\documentclass[11pt, reqno]{amsart}
\usepackage{amsmath, amsthm, amscd, amsfonts, amssymb, graphicx, latexsym, enumerate, color}
\usepackage[bookmarksnumbered, colorlinks, plainpages]{hyperref}
\newtheorem{theorem}{Theorem}[section]
\newtheorem{lemma}[theorem]{Lemma}
\newtheorem{proposition}[theorem]{Proposition}
\newtheorem{corollary}[theorem]{Corollary}
\theoremstyle{definition}
\newtheorem{definition}[theorem]{Definition}

\newcommand{\set}[1]{\left\{#1\right\}}

\numberwithin{equation}{section}
\makeatletter
\@namedef{subjclassname@2010}{%
\textup{2010} Mathematics Subject Classification}
\makeatother
\begin{document}
\setcounter{page}{1}
\title[(GORENSTEIN) INJECTIVE AND (GORENSTEIN) FLAT DIMENSIONS]{interplay between homological dimensions of a complex and its right derived section}
\author[C. JALALI]{CYRUS JALALI}
\address{C. Jalali, Faculty of mathematical sciences and computer, Kharazmi
University, 599 Taleghani Avenue, Tehran 15618, Iran}
\email{\textcolor[rgb]{0.00,0.00,0.84}{c.jalali@yahoo.com}}
\begin{abstract}
Let $(R,\mathfrak{m})$ be a commutative Noetherian local ring, $\mathfrak{a}$ be a proper ideal of $R$ and $M$ be an $R$-complex in $\mathrm{D}(R)$. We prove that if $M\in\mathrm{D}^f_\sqsubset(R)$ (respectively, $M\in\mathrm{D}^f_\sqsupset(R)$), then  $\mathrm{id}_R\mathbf{R}\Gamma_{\mathfrak{a}}(M)=\mathrm{id}_R M$ (respectively, $\mathrm{fd}_R\mathbf{R}\Gamma_{\mathfrak{a}}(M)=\mathrm{fd}_R M$). Next, it is proved that the right derived section functor of a complex $M\in\mathrm{D}_\sqsubset(R)$ ($R$ is not necessarily local) can be computed via a genuine left-bounded complex $G\simeq M$ of Gorenstein injective modules. We show that if $R$ has a dualizing complex and $M$ is an $R$-complex in $\mathrm{D}^f_\square(R)$, then $\mathrm{Gfd}_R\mathbf{R}\Gamma_{\mathfrak{a}}(M)=\mathrm{Gfd}_R M$ and $\mathrm{Gid}_R\mathbf{R}\Gamma_{\mathfrak{a}}(M)=\mathrm{Gid}_R M$. Also, we show that if $M$ is a relative Cohen-Macaulay $R$-module with respect to $\mathfrak{a}$ (respectively, Cohen-Macaulay $R$-module of dimension $n$), then $\mathrm{Gfd}_R\mathbf{H}^{\mathrm{ht_M\mathfrak{a}}}_{\mathfrak{a}}(M)=\mathrm{Gfd}_RM+n$ (respectively, $\mathrm{Gid}_R\mathbf{H}^n_{\mathfrak{m}}(M)=\mathrm{Gid}_RM-n$). The above results generalize some known results and provide characterizations of Gorenstein rings.
\end{abstract}
\maketitle
\section{Introduction}Throughout this paper, $R$ is a commutative Noetherian ring, $\mathfrak{a}$ is a proper ideal of $R$ and $M$ is an $R$-complex. The category of $R$-complexes is denoted $\mathrm{C}(R)$, and we use subscripts $\sqsubset$, $\sqsupset$ and $\square$ to denote genuine boundedness conditions. So, $\mathrm{C}_{\sqsupset}(R)$ is the full subcategory of $\mathrm{C}(R)$ of bounded below complexes (see {\cite[Definition 2.1.1]{C&Fo}}).
Also, the derived category is denoted $\mathrm{D}(R)$, and we use subscripts $\sqsubset$, $\sqsupset$ and $\square$ to denote homological boundedness conditions (see {\cite[Definition 4.1.15]{C&Fo}}).
The symbol $\simeq$ is the sign for isomorphism in $\mathrm{D}(R)$ and quasiisomorphisms in $\mathrm{C}(R)$. We also use superscript $f$ to signify that the homology modules are degreewise finitely generated. An $R$-complex $I$ is semiinjective if the functor $\mathrm{Hom}_{R}(-, I)$ converts injective quasiisomorphisms into surjective quasiisomorphisms. A semiinjective resolution of an $R$-complex $M$ is a semiinjective complex $I$ and a quasiisomorphism $M\overset{\simeq}{\longrightarrow}I$. For an $R$-complex $M$ the injective dimension $\mathrm{id}_R M$ is defined as
$$\mathrm{id}_R M=\mathrm{inf}\set{\ell\in \mathbb{Z}\mid \substack{\exists\ semiinjective\ R-complex\ I\ such\ that\\ M\simeq I\ in\ \mathrm{D}(R)\ and\ I_v=0\ for\ all\ v<-\ell}}.$$
Several of the main results of this paper involve the hypothesis that $R$ has a dualizing complex. A complex $D\in\mathrm{D}_\square^f(R)$ is dualizing for $R$ if it has finite injective dimension and the canonical morphism $\chi_M^R:R\rightarrow \mathbf{R}\mathrm{Hom}_{R}(D,D)$ is an isomorphism in $\mathrm{D}(R)$. If $R$ has a dualizing complex $D$, we may consider the functor $-^\dag=\mathbf{R}\mathrm{Hom}_{R}(-,D)$. The notion of Gorenstein injective module was introduced by E.E. Enochs and O.M.G. Jenda in {\cite{E&J}}. An $R$-module $M$ is said to be Gorenstein injective, if there exists a $\mathrm{Hom}_{R}(\mathcal{I},-)$ exact acyclic complex $E$ of injective $R$-modules such that $M=\mathrm{Ker}(E_0\rightarrow E_{-1})$. The Gorenstein injective dimension, $\mathrm{Gid}_RM$, of $M\in \mathrm{D}_\sqsubset(R)$ is defined to be the infimum of the set of integers $n$ such that there exists a complex $G\in\mathrm{C}_\sqsubset(R)$ consisting of Gorenstein injective modules satisfying $M\simeq G$ and $G_n=0$ for $n<-\ell$. Also, an R-complex $F$ is semiflat if the functor $-\otimes_RF$ preserves injective quasiisomorphisms. For an $R$-complex $M$ the flat dimension $\mathrm{fd}_R M$ is defined as
$$\mathrm{fd}_R M=\mathrm{inf}\set{n\in \mathbb{Z}\mid \substack{\exists\ semiflat\ R-complex\ F\ such\ that\\ F\simeq M\ in\ \mathrm{D}(R)\ and\ F_v=0\ for\ all\ v>n}}.$$
An $R$-module $M$ is said to be Gorenstein flat, if there exists an $\mathcal{I}\otimes_R-$ exact acyclic complex $F$ of of flat $R$-modules such that
$M=\mathrm{Ker}(F_0\rightarrow F_{-1})$. The Gorenstein flat dimension, $\mathrm{Gfd}_RM$, of $M\in \mathrm{D}_\sqsupset(R)$ is defined to be the infimum of the set of integers $n$ such that there exists a complex $F\in\mathrm{C}_\sqsupset(R)$ consisting of Gorenstein flat modules satisfying $M\simeq F$ and $G_n=0$ for $n<\ell$. Let $M$ be an $R$-complex in $\mathrm{D}(R)$. The right derived section functor of the complex $M$ is defined as $\mathbf{R}\Gamma_{\mathfrak{a}}(M)=\Gamma_{\mathfrak{a}}(E)$, where $E$ is a semiinjective resolution of $M$ (see {\cite{G&M}} and {\cite{Sp}}). If $\underline{x}= x_1,\cdots, x_r$ is a generating set for the ideal $\mathfrak{a}$ and $\check{C}_{\underline{x}}$ the corresponding $\check{\mathrm{C}}$ech complex, then $\mathbf{R}\Gamma_{\mathfrak{a}}(M)\simeq M\otimes_R^{\mathbf{L}}\check{C}_{\underline{x}}$ (see {\cite[Theorem 1.1(iv)]{Sc}}).\\
\indent It has been shown in {\cite[Theorem 6.5]{F}} that the right derived section functor (with support in any ideal $\mathfrak{a}$) sends complexes of finite flat dimension (respectively, finite injective dimension) to complexes of finite flat dimension (respectively, finite injective dimension). In section 2 we prove that if $(R,\mathfrak{m})$ is a local ring and $M\in\mathrm{D}^f_\sqsubset(R)$, then  $\mathrm{id}_R\mathbf{R}\Gamma_{\mathfrak{a}}(M)=\mathrm{id}_R M$ (see Theorem \ref{2th6}). It shows that the following statements are equivalent:
\begin{enumerate}[\upshape (i)]
\item $R$ is Gorenstein;
\item $\mathrm{id}_R\mathbf{R}\Gamma_{\mathfrak{a}}(R)=\mathrm{dim}(R)$ for any ideal $\mathfrak{a}$ of $R$;
\item $\mathrm{id}_R\mathbf{R}\Gamma_{\mathfrak{a}}(R)<\infty$ for some ideal $\mathfrak{a}$ of $R$.
\end{enumerate}
This provides a characterization of Gorenstein rings, which recovers {\cite[Corollary 2.7]{Z&Z}}.
Next, in \ref{2th8}, we prove that if $(R,\mathfrak{m})$ is a local ring and $M$ is an $R$-complex in $\mathrm{D}^f_\sqsubset(R)$ with $\mathrm{amp}\mathbf{R}\Gamma_{\mathfrak{a}}(M)=0$, then $\mathrm{id}_R\mathbf{H}^{-\mathrm{inf}\mathbf{R}\Gamma_{\mathfrak{a}}(M)}_{\mathfrak{a}}(M)=\mathrm{id}_R M+\mathrm{inf}\mathbf{R}\Gamma_{\mathfrak{a}}(M).$ Notice that, this result is a generalization of {\cite[Theorem 2.5]{Z&Z}}. Also, a flat version of \ref{2th6} is demonstrated. Indeed, it is shown, in Theorem \ref{3th6}, that if $(R,\mathfrak{m})$ is a local ring and $M\in\mathrm{D}^f_\sqsupset(R)$, then $\mathrm{fd}_R\mathbf{R}\Gamma_{\mathfrak{a}}(M)=\mathrm{fd}_R M$.\\
\indent It has been proved in {\cite[Theorem 5.9]{C&Fr&H}} that if $M\in \mathrm{D}_\square(R)$, then $$\mathrm{Gfd}_RM<\infty \Rightarrow \mathrm{Gfd}_R\mathbf{R}\Gamma_{\mathfrak{a}}(M)<\infty.$$
Moreover, if $R$ has a dualizing complex, the above implication may be reversed if $\mathfrak{a}$ is in the Jacobson radical of $R$ and $M\in \mathrm{D}^f_\square(R)$. We show that if $(R,\mathfrak{m})$ is a local ring and $M\in \mathrm{D}^f_\square(R)$, then
$$\mathrm{Gfd}_RM<\infty \Rightarrow \mathrm{Gfd}_R\mathbf{R}\Gamma_{\mathfrak{a}}(M)=\mathrm{Gfd}_RM,\ and$$
$$R\ admits\ a\ dualizing\ complex\ \Rightarrow \mathrm{Gfd}_R\mathbf{R}\Gamma_{\mathfrak{a}}(M)=\mathrm{Gfd}_R M.$$

Then, as a corollary, we prove that if $M$ is a relative Cohen-Macaulay $R$-module with respect to $\mathfrak{a}$, where is defined as in {\cite{Z&Z}}, and that $n=\mathrm{grade}(\mathfrak{a},M)$, then $\mathrm{Gfd}_R\mathbf{H}^n_{\mathfrak{a}}(M)=\mathrm{Gfd}_RM+n$.\\
\indent In section 3, first we prove that if $C$ is an $\Gamma_{\mathfrak{a}}$-acyclic $R$-complex in $\mathrm{C}_\sqsubset(R)$, then $\mathbf{R}\Gamma_{\mathfrak{a}}(C)\simeq \Gamma_{\mathfrak{a}}(C)$ (see Theorem \ref{4th5}).
It implies that the right derived section functor of a complex $M\in\mathrm{D}_\sqsubset(R)$ can be computed via a genuine left-bounded complex $G\simeq M$ of Gorenstein injective modules.\\
\indent Also, as a main result, we show that if $(R,\mathfrak{m})$ is a local ring admitting a dualizing complex and $M$ is an $R$-complex in $\mathrm{D}^f_\square(R)$, then $\mathrm{Gid}_R\mathbf{R}\Gamma_{\mathfrak{a}}(M)=\mathrm{Gid}_R M$ (see Theorem \ref{4th11}). It shows that the following statements are equivalent:
\begin{enumerate}[\upshape (i)]
\item $R$ is Gorenstein;
\item $\mathrm{Gid}_R\mathbf{R}\Gamma_{\mathfrak{a}}(R)=\mathrm{dim}(R)$ for any ideal $\mathfrak{a}$ of $R$;
\item $\mathrm{Gid}_R\mathbf{R}\Gamma_{\mathfrak{a}}(R)<\infty$ for some ideal $\mathfrak{a}$ of $R$.
\end{enumerate}
This provides a characterization of Gorenstein rings, which improves {\cite[Corollary 3.10]{Z&Z}} and {\cite[Theorem 2.6]{Yo}}, that is, we may prove them without assuming that $R$ is Cohen-Macaulay. Next, in Theorem \ref{4th12}, we prove a complex version of \ref{2th6}, which improves {\cite[Theorem 3.8]{Z&Z}}. As a corollary, in \ref{4co14}, we deduce that $\mathrm{Gid}_R\mathbf{H}^n_{\mathfrak{m}}(M)=\mathrm{Gid}_RM-n$, wherever $(R,\mathfrak{m})$ is a local ring and $M$ is a Cohen-Macaulay $R$-module with $\mathrm{dim}_R M=n$.
\section{right derived section functor, injective dimension and (Gorenstein) flat dimension}

The following lemma, which is an immediate consequence of {\cite[Corollary 3.4.4]{L}} and {\cite[Proposition 3.2.2]{L}}, determines the \textit{i}-th Bass number $\mu_{R_\mathfrak{p}}^i((\mathbf{R}\Gamma_{\mathfrak{a}}(M))_\mathfrak{p})$ of $(\mathbf{R}\Gamma_{\mathfrak{a}}(M))_\mathfrak{p}$ (see {\cite[Definition 6.1.18]{C&Fo}}).
\begin{lemma}\label{2co4}Let $(R,\mathfrak{m},k)$ be a local ring, and let $M$ be an $R$-complex in $\mathrm{D}(R)$. Then
$\mu_R^i(\mathbf{R}\Gamma_{\mathfrak{a}}(M))= \mu_R^i(M)$ for all $i\in \mathbb{Z}$; In particular, for every $\mathfrak{p}\in\mathbf{V}(\mathfrak{a})$ there is an equality $\mu_{R_\mathfrak{p}}^i((\mathbf{R}\Gamma_{\mathfrak{a}}(M))_\mathfrak{p})=\mu_{R_\mathfrak{p}}^i(M_\mathfrak{p})$ for all $i\in \mathbb{Z}$.
\end{lemma}
The following theorem, which is one of the main results of this section, provides a comparison between the injective dimensions of a complex and its right derived section functor.
\begin{theorem}\label{2th6}Let $(R,\mathfrak{m},k)$ be a local ring, and let $M$ be an $R$-complex in $\mathrm{D}^f_\sqsubset(R)$. Then $\mathrm{id}_R\mathbf{R}\Gamma_{\mathfrak{a}}(M)=\mathrm{id}_R M$.
\end{theorem}
\begin{proof} Let $s:=\mathrm{id}_R M<\infty$. Then, in view of {\cite[Lemma 6.1.19]{C&Fo}} and Lemma \ref{2co4}, $\mu^{i+s}_{R_\mathfrak{p}}((\mathbf{R}\Gamma_{\mathfrak{a}}(M))_\mathfrak{p})=0$ for all $\mathfrak{p}\in\mathrm{Spec}(R)$ and for all $i>0$. Therefore, it follows from {\cite[Lemma 6.1.19]{C&Fo}} that $\mathrm{id}_R\mathbf{R}\Gamma_{\mathfrak{a}}(M)\leq s$.\\
\indent For the opposite inequality, let $t:=\mathrm{id}_R\mathbf{R}\Gamma_{\mathfrak{a}}(M)<\infty $. Then, by {\cite[Theorem 5.1.6]{C&Fo}}, $\mathrm{inf}\mathbf{R}\mathrm{Hom}_{R}(T,\mathbf{R}\Gamma_{\mathfrak{a}}(M))\geq -t$ for all cyclic R-modules $T$.  Hence, in view of {\cite[Proposition 3.2.2]{L}}, there are isomorphisms
$$\mathbf{H}_{{-t-i}}(\mathbf{R}\mathrm{Hom}_{R}(k,M))\cong \mathbf{H}_{-t-i}(\mathbf{R}\mathrm{Hom}_{R}(k,\mathbf{R}\Gamma_{\mathfrak{a}}(M)))\cong 0$$
for all $i>0$. Therefore $-\mathrm{inf}\mathbf{R}\mathrm{Hom}_{R}(k,M)\leq t$, and so $\mathrm{id}_R M\leq t$ by {\cite[Theorem 6.1.13]{C&Fo}}.
\end{proof}
The following corollary, which recovers {\cite[Corollary 2.7]{Z&Z}}, is an immediate consequence of the previous Theorem.
\begin{corollary}\label{2co7}Let $(R,\mathfrak{m})$ be a local ring. Then the following statements are equivalent:
\begin{enumerate}[\upshape (i)]
\item $R$ is Gorenstein;
\item $\mathrm{id}_R\mathbf{R}\Gamma_{\mathfrak{a}}(R)=\mathrm{dim}(R)$ for any ideal $\mathfrak{a}$ of $R$;
\item $\mathrm{id}_R\mathbf{R}\Gamma_{\mathfrak{a}}(R)<\infty$ for some ideal $\mathfrak{a}$ of $R$.
\end{enumerate}
\end{corollary}
The following theorem,  which is an immediate consequence of Theorem \ref{2th6}, is a generalization of {\cite[Theorem 2.5]{Z&Z}}.
\begin{theorem}\label{2th8}Let $(R,\mathfrak{m})$ be a local ring. Suppose that $M$ is an $R$-complex in $\mathrm{D}_\square^f(R)$ such that $\mathrm{amp}\mathbf{R}\Gamma_{\mathfrak{a}}(M)=0$. Then $$\mathrm{id}_R\mathbf{H}^{-\mathrm{inf}\mathbf{R}\Gamma_{\mathfrak{a}}(M)}_{\mathfrak{a}}(M)=\mathrm{id}_R M+\mathrm{inf}\mathbf{R}\Gamma_{\mathfrak{a}}(M).$$
\end{theorem}
\begin{proof}Let $n:=-\mathrm{inf}\mathbf{R}\Gamma_{\mathfrak{a}}(M)$. Since $\mathbf{R}\Gamma_{\mathfrak{a}}(M)\simeq\mathbf{H}_{-n}(\mathbf{R}\Gamma_{\mathfrak{a}}(M))$, there is an equality
\begin{flalign*}
\mathrm{id}_R\mathbf{R}\Gamma_{\mathfrak{a}}(M)
&=\mathrm{id}_R\Sigma^n\mathbf{H}_{-n}(\mathbf{R}\Gamma_{\mathfrak{a}}(M))+n.
\end{flalign*}
But $\Sigma^n\mathbf{R}\Gamma_{\mathfrak{a}}(M)$ is equivalent to the module $\Sigma^n\mathbf{H}_{-n}(\mathbf{R}\Gamma_{\mathfrak{a}}(M))$ in the category of R-modules. So, we may identify $\Sigma^n\mathbf{H}_{-n}(\mathbf{R}\Gamma_{\mathfrak{a}}(M))$ with $\mathbf{H}^n_{\mathfrak{a}}(M)$. Hence $\mathrm{id}_R\mathbf{H}^n_{\mathfrak{a}}(M)=\mathrm{id}_R\mathbf{R}\Gamma_{\mathfrak{a}}(M)-n$. The desired
equality now follows from Theorem \ref{2th6}.
\end{proof}
In the following we use the notion of a semifree resolution. A semifree resolution of an $R$-complex $M$ is a semifree complex $F$ (see {\cite[Definition 3.1.1]{C&Fo}}) and a quasiisomorphism $F\overset{\simeq}{\longrightarrow}M$.
\begin{lemma}\label{3pr1}Let $(R,\mathfrak{m},k)$ be a local ring, and let $M$ be an $R$-complex in $\mathrm{D}(R)$. Then $$k\otimes_R^{\mathbf{L}}\mathbf{R}\Gamma_{\mathfrak{a}}(M)\simeq k\otimes_R^{\mathbf{L}}M.$$
\end{lemma}
\begin{proof}Let $F$ be a semifree resolution of the residue field k, and let $\underline{x}= x_1,\cdots, x_r$ be a generating set for the ideal $\mathfrak{a}$ and $\check{C}_{\underline{x}}$ be the $\check{\mathrm{C}}$ech complex with respect to $\underline{x}$. But, as in the proof of {\cite[Lemma 2.4]{M&Sc}}, there exists a quasiisomorphism $F\otimes_R\check{C}_{\underline{x}}\simeq F$ in $\mathrm{C}(R)$. The result now follows, since $k\otimes_R^{\mathbf{L}}\check{C}_{\underline{x}}\simeq k$.
\end{proof}
The following lemma, which is an immediate consequence of {\cite[Corollary 3.4.4]{L}} and \ref{3pr1}, determines the \textit{i}-th Betti number $\beta_i^{R_\mathfrak{p}}((\mathbf{R}\Gamma_{\mathfrak{a}}(M))_\mathfrak{p})$ of $(\mathbf{R}\Gamma_{\mathfrak{a}}(M))_\mathfrak{p}$ (see {\cite[Definition 6.1.14]{C&Fo}}).
\begin{lemma}\label{3co2}Let $(R,\mathfrak{m},k)$ be a local ring, and let $M$ be an $R$-complex in $\mathrm{D}(R)$. Then
$\beta_i^R(\mathbf{R}\Gamma_{\mathfrak{a}}(M))= \beta_i^R(M)$ for all $i\in \mathbb{Z}$; In particular, for every $\mathfrak{p}\in\mathbf{V}(\mathfrak{a})$ there is an equality $\beta_i^{R_\mathfrak{p}}((\mathbf{R}\Gamma_{\mathfrak{a}}(M))_\mathfrak{p})=\beta_i^{R_\mathfrak{p}}(M_\mathfrak{p})$
for all $i\in \mathbb{Z}$.
\end{lemma}
\begin{theorem}\label{3th6}Let $(R,\mathfrak{m},k)$ be a local ring, and let $\mathfrak{a}$ be a proper ideal of $R$. Suppose that $M$ is an $R$-complex in $\mathrm{D}^f_\sqsupset(R)$. Then $\mathrm{fd}_R\mathbf{R}\Gamma_{\mathfrak{a}}(M)=\mathrm{fd}_R M$.
\end{theorem}
\begin{proof} Let $s:=\mathrm{fd}_R M<\infty$. Then, in view of {\cite[Lemma 6.1.15]{C&Fo}} and Lemma \ref{3co2}, $\beta_{i+s}^{R_\mathfrak{p}}((\mathbf{R}\Gamma_{\mathfrak{a}}(M))_\mathfrak{p})=0$ for all $\mathfrak{p}\in\mathrm{Spec}(R)$ and for all $i>0$. Therefore, it follows from {\cite[Lemma 6.1.15]{C&Fo}} that $\mathrm{fd}_R\mathbf{R}\Gamma_{\mathfrak{a}}(M)\leq s$.\\
\indent For the opposite inequality, let $t:=\mathrm{fd}_R\mathbf{R}\Gamma_{\mathfrak{a}}(M)<\infty $. Then, by {\cite[Theorem 5.1.9]{C&Fo}}, $\mathrm{sup}T\otimes_R^{\mathbf{L}}\mathbf{R}\Gamma_{\mathfrak{a}}(M))\geq t$ for all cyclic R-modules $T$.  Hence, in view of Lemma \ref{3pr1}, there are isomorphisms
$$\mathbf{H}_{{t+i}}(k\otimes_R^{\mathbf{L}}M)\cong \mathbf{H}_{t+i}(k\otimes_R^{\mathbf{L}}\mathbf{R}\Gamma_{\mathfrak{a}}(M))\cong 0$$
for all $i>0$. Therefore $\mathrm{sup}k\otimes_R^{\mathbf{L}}M\leq t$, and so $\mathrm{fd}_R M\leq t$ by {\cite[Theorem 5.2.13]{C&Fo}}.
\end{proof}
The following theorem, which is an immediate consequence of Theorem \ref{3th6}, is a flat version of \ref{2th8}.
\begin{theorem}\label{3th7}Let $(R,\mathfrak{m})$ be a local ring. Suppose that $M$ is an $R$-complex in $\mathrm{D}_\square^f(R)$ such that $\mathrm{amp}\mathbf{R}\Gamma_{\mathfrak{a}}(M)=0$. Then $$\mathrm{fd}_R\mathbf{H}^{-\mathrm{inf}\mathbf{R}\Gamma_{\mathfrak{a}}(M)}_{\mathfrak{a}}(M)=\mathrm{fd}_R M-\mathrm{inf}\mathbf{R}\Gamma_{\mathfrak{a}}(M).$$
\end{theorem}
\begin{proof} Straightforward verification similar to the proof of Theorem \ref{2th8}.
\end{proof}
In the rest of this section, we make a comparison between the Gorenstein flat dimensions of a complex and its right derived section functor.
\begin{proposition}\label{3pr8}Suppose that $M$ is an $R$-complex in $\mathrm{D}_\square(R)$. Then $$\mathrm{Gfd}_R\mathbf{R}\Gamma_{\mathfrak{a}}(M)\leq\mathrm{Gfd}_R M.$$
\end{proposition}
\begin{proof}Notice that if $\mathrm{Gfd}_RM=\infty$, then there is nothing to prove. So, we may assume that $\mathrm{Gfd}_RM<\infty$. Hence, it follows from {\cite[Theorem 5.9]{C&Fr&H}} that $\mathrm{Gfd}_R\mathbf{R}\Gamma_{\mathfrak{a}}(M)<\infty$. Now, by {\cite[Theorem 8.8]{I&S}}, there exists  $\mathfrak{p}\in\mathbf{V}(\mathfrak{a})$ such that $$\mathrm{Gfd}_R\mathbf{R}\Gamma_{\mathfrak{a}}(M)=\mathrm{depth}R_\mathfrak{p}-\mathrm{depth}_{R_\mathfrak{p}}(\mathbf{R}\Gamma_{\mathfrak{a}}(M))_\mathfrak{p}.$$
But $\mathrm{depth}_{R_\mathfrak{p}}(\mathbf{R}\Gamma_{\mathfrak{a}}(M))_\mathfrak{p}=\mathrm{depth}_{R_\mathfrak{p}}M_\mathfrak{p}$. It follows, again by {\cite[Theorem 8.8]{I&S}}, that  $$\mathrm{Gfd}_R\mathbf{R}\Gamma_{\mathfrak{a}}(M)=\mathrm{depth}R_\mathfrak{p}-\mathrm{depth}_{R_\mathfrak{p}}M_\mathfrak{p}\leq\mathrm{Gfd}_RM$$ as desired.
\end{proof}
\begin{proposition}\label{3pr9}Let $(R,\mathfrak{m},k)$ be a local ring, and let $M$ be an $R$-complex in $\mathrm{D}_\square^f(R)$ such that $\mathrm{Gfd}_RM<\infty$. Then $$\mathrm{Gfd}_R M\leq\mathrm{Gfd}_R\mathbf{R}\Gamma_{\mathfrak{a}}(M).$$
\end{proposition}
\begin{proof}By {\cite[Theorem 5.9]{C&Fr&H}}, $\mathrm{Gfd}_R\mathbf{R}\Gamma_{\mathfrak{a}}(M)<\infty$. Hence, by {\cite[Theorem 8.7]{I&S}}, there are equalities
\begin{flalign*}
\mathrm{sup}(E(k)\otimes_R^{\mathbf{L}}\mathbf{R}\Gamma_{\mathfrak{a}}(M))
&=\mathrm{depth}R-\mathrm{depth}_R\mathbf{R}\Gamma_{\mathfrak{a}}(M)\\
&=\mathrm{depth}R-\mathrm{depth}_RM=\mathrm{sup}(E(k)\otimes_R^{\mathbf{L}}M).
\end{flalign*}
Since $M\in \mathrm{D}_\square^f(R)$, $\mathrm{sup}(E(k)\otimes_R^{\mathbf{L}}M)=\mathrm{Gfd}_R M$. The result now follows from the fact that
$\mathrm{Gfd}_R\mathbf{R}\Gamma_{\mathfrak{a}}(M)\geq \mathrm{sup}(E(k)\otimes_R^{\mathbf{L}}\mathbf{R}\Gamma_{\mathfrak{a}}(M))$ (see {\cite[Corollary 3.6]{C&Fr&H}}).
\end{proof}
The following theorem is a Gorenstein flat version of Theorem \ref{3th6}.
\begin{theorem}\label{3th10}Let $(R,\mathfrak{m})$ be a local ring, and let $M$ be an $R$-complex in $\mathrm{D}_\square^f(R)$.
\begin{enumerate}[\upshape (i)]
\item If $\mathrm{Gfd}_R M<\infty$, then $\mathrm{Gfd}_R\mathbf{R}\Gamma_{\mathfrak{a}}(M)=\mathrm{Gfd}_R M$.
\item If $R$ admits a dualizing complex, then $\mathrm{Gfd}_R\mathbf{R}\Gamma_{\mathfrak{a}}(M)=\mathrm{Gfd}_R M$.
\end{enumerate}
\end{theorem}
\begin{proof}(i) A straightforward application of Proposition \ref{3pr8} and Proposition \ref{3pr9}.\\
\indent (ii) In view of part (i), we may assume that $\mathrm{Gfd}_R\mathbf{R}\Gamma_{\mathfrak{a}}(M)<\infty$. Hence, by {\cite[Theorem 5.9]{C&Fr&H}}, $\mathrm{Gfd}_R M<\infty$. The desired equality now follows from Proposition \ref{3pr8} and Proposition \ref{3pr9}.
\end{proof}
\begin{corollary}\label{3co12}Let $(R,\mathfrak{m})$ be a local ring. Suppose that $M$ is relative Cohen-Macaulay with respect to $\mathfrak{a}$ and that
$n=\mathrm{grade}(\mathfrak{a},M)$. Then $\mathrm{Gfd}_R\mathbf{H}^n_{\mathfrak{a}}(M)=\mathrm{Gfd}_RM+n$.
\end{corollary}
\begin{proof} Notice that $(\widehat{R},\widehat{\mathfrak{m}})$ is a local ring admitting a dualizing complex and $M\otimes_R \widehat{R}$ is a relative Cohen-Macaulay $\widehat{R}$-module with respect to $\mathfrak{a}\widehat{R}$ and that
$\mathrm{grade}(\mathfrak{a}\widehat{R},M\otimes_R \widehat{R})=n$. Hence, in view of {\cite[Theorem 4.27]{C&Fo&H}}, we may assume that $R$ is complete; and so it has a dualizing complex. The result therefore follows from Theorem \ref{3th10}.
\end{proof}
\section{right derived section functor and Gorenstein injective dimension}
In this section the category of $R$-modules is denoted $\mathcal{C}(R)$. Recall from {\cite[Exercise 4.1.2]{B&S}} that the local cohomology modules of $R$-module $M$ with respect to $\mathfrak{a}$ can be calculated by an $\Gamma_{\mathfrak{a}}$-acyclic resolution of $M$. First, we prove the complex version of it.
\begin{definition}\label{4de1}(see {\cite[5.7.9]{W}}) Let $\mathrm{F}:\mathcal{C}(R)\rightarrow\mathcal{C}(R)$ be a left exact functor, and assume that $M$ is an $R$-complex in $\mathrm{C}_\sqsubset(R)$. If $0\rightarrow M\rightarrow C_{*,0}\rightarrow C_{*,1}\rightarrow\cdots\rightarrow C_{*,q}\rightarrow$ is a Cartan-Eilenberg injective resolution of $M$, where is defined as in {\cite[\S10.5]{R}}, define $\mathbb{R}^i(\mathrm{F}M)$ to be $\mathbf{H}_i(\mathrm{Tot}(\mathrm{F}C))$.
\end{definition}
\begin{lemma}\label{4pr2}Let $M$ and $\acute{M}$ be two $R$-complexes in $\mathrm{C}_\sqsubset(R)$, and let $\zeta:M\rightarrow \acute{M}$ be a morphism of $R$-complexes. Suppose that $0\rightarrow M\rightarrow C_{*,0}\rightarrow C_{*,1}\rightarrow\cdots\rightarrow C_{*,q}\rightarrow$ and $0\rightarrow \acute{M}\rightarrow \acute{C}_{*,0}\rightarrow \acute{C}_{*,1}\rightarrow\cdots\rightarrow \acute{C}_{*,q}\rightarrow$ are Cartan-Eilenberg injective resolutions of $M$ and $\acute{M}$, respectively. Then there exists a sequence $\{\zeta_{*,q}\}_{q\in\mathbb{N}_0}$ of morphisms $\zeta_{*,q}:C_{*,q}\rightarrow \acute{C}_{*,q}$ of $R$-complexes over $\zeta$.
\end{lemma}
\begin{proof}A straightforward application of {\cite[Theorem 19]{N}}.
\end{proof}
\begin{lemma}\label{4le3}Let $\mathrm{F}:\mathcal{C}(R)\rightarrow\mathcal{C}(R)$ be a left exact functor, and let $M$ and $\acute{M}$ be two $R$-complexes. Then
\begin{enumerate}[\upshape (i)]
\item Any quasiisomorphism $\zeta:M\rightarrow \acute{M}$ induces isomorphism $$\mathbb{R}^i(\mathrm{F}M)\cong \mathbb{R}^i(\mathrm{F}\acute{M})$$
for all $i\in \mathbb{Z}$; and
\item If $M$ is $\mathrm{F}$-acyclic $R$-complex in $\mathrm{C}_\sqsubset(R)$, that is, $M_i$ is $\mathrm{F}$-acyclic for all $i\in \mathbb{Z}$, then $$\mathbb{R}^i(\mathrm{F}M)= \mathbf{H}_i(\mathrm{F}M)$$ for all $i\in \mathbb{Z}$.
\end{enumerate}
\end{lemma}
\begin{proof}
Straightforward verification similar to the proof of {\cite[Corollary 5.7.7]{W}}.
\end{proof}
The following theorem, which is one of the main results of this section, enables us to prove some interesting results.
\begin{theorem}\label{4th4}Let $\mathrm{F}:\mathcal{C}(R)\rightarrow\mathcal{C}(R)$ be a left exact functor, and let $M$ be an $\mathrm{F}$-acyclic $R$-complex in $\mathrm{C}_\sqsubset(R)$. Assume that $I$ is $\mathrm{F}$-acyclic and $\mathrm{F}(I)$ is injective for every injective $R$-module I. Then  $\mathrm{F}(M)\simeq \mathrm{F}(E)$, for every semiinjective resolution $E\in \mathrm{C}_\sqsubset(R)$ of $M$.
\end{theorem}
\begin{proof}Let $E$ be a semiinjective resolution of $M$ with $E_v=0$ for $v>\mathrm{sup}M$, and let $\zeta:M\overset{\simeq}{\longrightarrow}E$ be an quasiisomorphism. By {\cite[Theorem 10.45]{R}}, there exist Cartan-Eilenberg injective resolutions $0\rightarrow M\rightarrow C_{*,0}\rightarrow C_{*,1}\rightarrow\cdots\rightarrow C_{*,q}\rightarrow$ and $0\rightarrow E\rightarrow \acute{C}_{*,0}\rightarrow \acute{C}_{*,1}\rightarrow\cdots\rightarrow \acute{C}_{*,q}\rightarrow$. Hence, in view of Lemma \ref{4pr2}, there is a sequence $\{\zeta_{*,q}\}_{q\in\mathbb{N}_0}$ of morphisms of $R$-complexes such that the diagram
$$\begin{CD}
0 @>>> E @>>> \acute{C}_{*,0} @>>> \acute{C}_{*,1} @>>>\cdots@>>> \acute{C}_{*,q} @>>>\\
@. @AA{\zeta}A @AA{\zeta_{*,0}}A @AA{\zeta_{*,1}}A @. @AA{\zeta_{*,q}}A\\
0 @>>> M @>>> C_{*,0} @>>> C_{*,1} @>>>\cdots@>>> C_{*,q} @>>>\\
\end{CD}$$
commutes in $\mathrm{C}(R)$. By Lemma \ref{4le3}(ii), $\mathbb{R}^p(\mathrm{F}(M))= \mathbf{H}_p(\mathrm{F}(M))$
for all $p\in \mathbb{Z}$, so that the natural morphism $\mathrm{F}(M)\rightarrow \mathrm{Tot}(\mathrm{F}(C))$ is a quasiisomorphism. Similarly, $\mathbb{R}^p(\mathrm{F}(E))= \mathbf{H}_p(\mathrm{F}(E))$ for all $p\in \mathbb{Z}$, so that the natural morphism $\mathrm{F}(E)\rightarrow \mathrm{Tot}(\mathrm{F}(\acute{C}))$ is a quasiisomorphism. Thus, by {\cite[Proposition 3.3.5(a)]{C&Fo}}, there exists a quasiisomorphism $\mathrm{Tot}(\mathrm{F}(\acute{C}))\rightarrow \mathrm{F}(E)$, since $\mathrm{F}(E)$ is injective.\\
\indent But, by Lemma \ref{4le3}(i), there are isomorphisms $\mathbb{R}^p(\mathrm{F}(M))\cong \mathbb{R}^p(\mathrm{F}(E))$
for all $p\in \mathbb{Z}$. Hence the morphism $\zeta_*:\mathrm{Tot}(\mathrm{F}(C))\longrightarrow \mathrm{Tot}(\mathrm{F}(\acute{C}))$ is a quasiisomorphism, where $$\zeta_n=\sum_{p+q=n}\mathrm{F}(\zeta_{p,q}):\mathrm{Tot}(\mathrm{F}(C))_n\rightarrow \mathrm{Tot}(\mathrm{F}(\acute{C}))_n$$ for all $n\in \mathbb{Z}$. Therefore, there are quasiisomorphisms
$$\mathrm{F}(M)\overset{\simeq}{\longrightarrow}\mathrm{Tot}(\mathrm{F}(C))\overset{\simeq}{\longrightarrow}\mathrm{Tot}(\mathrm{F}(\acute{C}))\overset{\simeq}{\longrightarrow}\mathrm{F}(E).$$
Now $\mathrm{F}(M)\simeq \mathrm{F}(E)$ as desired.
\end{proof}
The next theorem, which offers an application of the previous theorem, is a complex version of {\cite[Exercise 4.1.2]{B&S}}.
\begin{theorem}\label{4th5}Let $C$ be an $\Gamma_{\mathfrak{a}}$-acyclic $R$-complex in $\mathrm{C}_\sqsubset(R)$. Then $\mathbf{R}\Gamma_{\mathfrak{a}}(C)\simeq \Gamma_{\mathfrak{a}}(C)$.
\end{theorem}
The next corollary shows that if $M\in\mathrm{D}_\sqsubset(R)$, then $\mathbf{R}\Gamma_{\mathfrak{a}}(M)$ can be computed via a genuine left-bounded complex $G\simeq M$ of Gorenstein injective modules. Also, notice that {\cite[Theorem 3.4]{S}} is an immediate consequence of this fact.
\begin{corollary}\label{4th7}Let $M$ be an $R$-complex in $\mathrm{D}_\sqsubset(R)$, and let $G\in \mathrm{C}_\sqsubset(R)$ be an $R$-complex of Gorenstein injective modules such that $M\simeq G$. Then $\mathbf{R}\Gamma_{\mathfrak{a}}(M)\simeq \Gamma_{\mathfrak{a}}(G)$.
\end{corollary}
\begin{proof}Since by {\cite[Lemma 1.1]{Y}} every Gorenstein injective module is $\Gamma_{\mathfrak{a}}$-acyclic, the result follows from Theorem \ref{4th5}.
\end{proof}
It has been proved in {\cite[Corollary 3.3]{S}} that if $R$ admits a dualizing complex and $M\in \mathrm{D}_\sqsubset(R)$, then $\mathrm{Gid}_R\mathbf{R}\Gamma_{\mathfrak{a}}(M)\leq\mathrm{Gid}_R M$. The following proposition together with {\cite[Theorem 5.9]{C&Fr&H}} recover this result.
\begin{proposition}\label{4pr8}Suppose that $M$ is an $R$-complex in $\mathrm{D}_\sqsubset(R)$ such that $\mathrm{Gid}_R\mathbf{R}\Gamma_{\mathfrak{a}}(M)<\infty$. Then $$\mathrm{Gid}_R\mathbf{R}\Gamma_{\mathfrak{a}}(M)\leq\mathrm{Gid}_R M.$$
\end{proposition}
\begin{proof}Notice that if $\mathrm{Gid}_RM=\infty$, then there is nothing to prove. So, we may assume that $\mathrm{Gid}_RM<\infty$. By {\cite[Theorem 2.2]{C&S}}, there exists  $\mathfrak{p}\in\mathbf{V}(\mathfrak{a})$ such that $$\mathrm{Gid}_R\mathbf{R}\Gamma_{\mathfrak{a}}(M)=\mathrm{depth}R_\mathfrak{p}-\mathrm{width}_{R_\mathfrak{p}}(\mathbf{R}\Gamma_{\mathfrak{a}}(M))_\mathfrak{p}.$$
But $\mathrm{width}_{R_\mathfrak{p}}(\mathbf{R}\Gamma_{\mathfrak{a}}(M))_\mathfrak{p}=\mathrm{width}_{R_\mathfrak{p}}M_\mathfrak{p}$. It follows, again by {\cite[Theorem 2.2]{C&S}}, that $$\mathrm{Gid}_R\mathbf{R}\Gamma_{\mathfrak{a}}(M)=\mathrm{depth}R_\mathfrak{p}-\mathrm{width}_{R_\mathfrak{p}}M_\mathfrak{p}\leq\mathrm{Gid}_RM$$ as desired.
\end{proof}
\begin{proposition}\label{4pr9}Let $(R,\mathfrak{m})$ be a local ring, and let $M$ be an $R$-complex in $\mathrm{D}_\square^f(R)$ such that $\mathrm{Gid}_RM<\infty$. Then $$\mathrm{Gid}_R M\leq\mathrm{Gid}_R\mathbf{R}\Gamma_{\mathfrak{a}}(M).$$
\end{proposition}
\begin{proof}By {\cite[Proposition 6.3]{C&Fr&H}}, there is an inequality
$$\mathrm{Gid}_R\mathbf{R}\Gamma_{\mathfrak{a}}(M)\geq \mathrm{depth}R-\mathrm{width}_R\mathbf{R}\Gamma_{\mathfrak{a}}(M).$$
But $\mathrm{width}_R\mathbf{R}\Gamma_{\mathfrak{a}}(M)=\mathrm{width}_RM$. Thus, we have
\begin{flalign*}
\mathrm{Gid}_R\mathbf{R}\Gamma_{\mathfrak{a}}(M)
&\geq\mathrm{depth}R-\mathrm{width}_RM\\
&=\mathrm{depth}R-\mathrm{inf}\ M.
\end{flalign*}
The result therefore follows from {\cite[Corollary 2.3]{C&S}}.
\end{proof}
The following theorem, which is a Gorenstein injective version of Theorem \ref{2th6}, is one of the main results of this section.
\begin{theorem}\label{4th11}Let $(R,\mathfrak{m})$ be a local ring admitting a dualizing complex, and let $M$ be an $R$-complex in $\mathrm{D}_\square^f(R)$. Then $\mathrm{Gid}_R\mathbf{R}\Gamma_{\mathfrak{a}}(M)=\mathrm{Gid}_R M$.
\end{theorem}
\begin{proof} A straightforward application of {\cite[Theorem 5.9]{C&Fr&H}}, Proposition \ref{4pr8} and Proposition \ref{4pr9}.
\end{proof}
The next theorem, which is a Gorenstein injective version of Theorem \ref{2th8}, recovers {\cite[Theorem 3.8]{Z&Z}}.
\begin{theorem}\label{4th12}Let $(R,\mathfrak{m})$ be a local ring admitting a dualizing complex. Suppose that $M$ is an $R$-complex in $\mathrm{D}_\square^f(R)$ such that $\mathrm{amp}\mathbf{R}\Gamma_{\mathfrak{a}}(M)=0$. Then $$\mathrm{Gid}_R\mathbf{H}^{-\mathrm{inf}\mathbf{R}\Gamma_{\mathfrak{a}}(M)}_{\mathfrak{a}}(M)=\mathrm{Gid}_R M+\mathrm{inf}\mathbf{R}\Gamma_{\mathfrak{a}}(M).$$
\end{theorem}
\begin{proof} Follows from Theorem \ref{4th11} (similar to the proof of Theorem \ref{2th8}).
\end{proof}
The following corollary, which improves {\cite[Corollary 3.9]{Z&Z}}, is a consequence of the previous Theorem.
\begin{corollary}\label{4co14}Let $(R,\mathfrak{m})$ be a local ring, and let $M$ be a Cohen-Macaulay $R$-module with $\mathrm{dim}_R M=n$. Then the following statements hold.
\begin{enumerate}[\upshape (i)]
\item $\mathrm{Gid}_{\widehat{R}}\mathbf{R}\Gamma_{\mathfrak{a}\widehat{R}}(M\otimes_R \widehat{R})=\mathrm{Gid}_RM$.
\item $\mathrm{Gid}_R\mathbf{H}^n_{\mathfrak{m}}(M)=\mathrm{Gid}_RM-n$.
\end{enumerate}
\end{corollary}
\begin{proof} Notice that $(\widehat{R},\widehat{\mathfrak{m}})$ is a local ring admitting a dualizing complex and $M\otimes_R \widehat{R}$ is a Cohen-Macaulay $\widehat{R}$-module of dimension $n$.\\
\indent (i) A straightforward application of Theorem \ref{4th11} and {\cite[Theorem 3.24]{C&Fo&H}}.\\
\indent (ii) There is an inequality
$$\mathrm{Gid}_{\widehat{R}}\mathbf{R}\Gamma_{\widehat{\mathfrak{m}}}(M\otimes_R \widehat{R})=\mathrm{Gid}_{\widehat{R}}\mathbf{H}^n_{\widehat{\mathfrak{m}}}(M\otimes_R \widehat{R})+n.$$
But, in view of {\cite[Lemma 3.6]{S}}, $\mathrm{Gid}_{\widehat{R}}\mathbf{H}^n_{\widehat{\mathfrak{m}}}(M\otimes_R \widehat{R})=\mathrm{Gid}_R\mathbf{H}^n_{\mathfrak{m}}(M)$. The result now follows from part (i).
\end{proof}
The next corollary provides a characterization of Gorenstein rings, which together with Corollary \ref{2co7} show that {\cite[Theorem 2.6]{Yo}} and {\cite[Corollary 3.10]{Z&Z}} hold without assuming that $R$ is Cohen-Macaulay.
\begin{corollary}\label{4co13}Let $(R,\mathfrak{m})$ be a local ring admitting a dualizing complex. Then the following statements are equivalent:
\begin{enumerate}[\upshape (i)]
\item $R$ is Gorenstein;
\item $\mathrm{Gid}_R\mathbf{R}\Gamma_{\mathfrak{a}}(R)=\mathrm{dim}(R)$ for any ideal $\mathfrak{a}$ of $R$;
\item $\mathrm{Gid}_R\mathbf{R}\Gamma_{\mathfrak{a}}(R)<\infty$ for some ideal $\mathfrak{a}$ of $R$.
\end{enumerate}
\end{corollary}
\begin{proof}A straightforward application of Theorem \ref{4th11} and {\cite[Proposition 3.11]{C&Fo&H}}.
\end{proof}

\end{document}